\documentclass[12pt,twoside,reqno]{amsart}

\usepackage{version}         

 \usepackage{color}

\usepackage{multirow}
\usepackage[OT1]{fontenc}
\usepackage{type1cm}
\usepackage{amssymb}
\usepackage{mathrsfs}

\usepackage{multicol}
\usepackage{multirow}
\usepackage[english]{babel}
\usepackage[latin1]{inputenc}
\usepackage[T1]{fontenc}
\usepackage{palatino}
\usepackage{amsfonts}
\usepackage{amsmath}
\usepackage{amssymb}
\usepackage{amsthm}
\usepackage{bbm}
\usepackage{extarrows}
\usepackage{graphicx}
\usepackage[usenames,dvipsnames]{xcolor}
\usepackage{wrapfig}
\usepackage{type1cm}
\usepackage[bf]{caption}
\usepackage{esint}
\usepackage{version}
\usepackage[colorlinks = true, citecolor = black, linkcolor = black, urlcolor = black, pdfstartview=FitH]{hyperref}
\usepackage{caption}
\usepackage{tikz}
\usepackage{bbding, wasysym}
\usepackage{enumitem}

\usepackage[left=2.3cm,top=3cm,right=2.3cm]{geometry}
\numberwithin{equation}{section}

\theoremstyle{plain}
\newtheorem{theorem}{Theorem}[section]

\newtheorem{corollary}[theorem]{Corollary}
\newtheorem{proposition}[theorem]{Proposition}

\theoremstyle{remark}
\newtheorem{remark}[theorem]{Remark}

\newtheorem{example}[theorem]{Example}

\newtheorem*{ack}{Acknowledgement}

\theoremstyle{definition}
\newtheorem{definition}[theorem]{Definition}

\newcommand{\C}{\mathbb{C}}

\newcommand{\N}{\mathbb{N}}

\newcommand{\cO}{\mathcal{O}}

\newcommand{\Z}{\mathbb{Z}}

\renewcommand{\epsilon}{\varepsilon}
\renewcommand{\rho}{\varrho}

\renewcommand{\k}{\mathtt{k}}

\begin{document}

\title{Infinite products related to generalized Thue-Morse sequences}

\author{Yao-Qiang Li}
\address{Institut de Math\'ematiques de Jussieu-Paris Rive Gauche \\
         Sorbonne Universit\'e \\
         Paris, 75005 \\
         France}
\email{yaoqiang.li@etu.upmc.fr\quad yaoqiang.li@imj-prg.fr}
\address{School of Mathematics \\
         South China University of Technology \\
         Guangzhou, 510641 \\
         P.R. China}
\email{scutyaoqiangli@qq.com\quad scutyaoqiangli@gmail.com}

\subjclass[2010]{Primary 11B83; Secondary 11B85, 68R15.}
\keywords{Thue-Morse sequence, morphic sequence, closed formulas for infinite products, Woods and Robbins product}
\date{}

\begin{abstract}
Given an integer $q\ge2$ and $\theta_1,\cdots,\theta_{q-1}\in\{0,1\}$, let $(\theta_n)_{n\ge0}$ be the generalized Thue-Morse sequence, defined to be the unique fixed point of the morphism
$$0\mapsto0\theta_1\cdots\theta_{q-1}$$
$$1\mapsto1\overline{\theta}_1\cdots\overline{\theta}_{q-1}$$
beginning with $\theta_0:=0$, where $\overline{0}:=1$ and $\overline{1}:=0$. For rational functions $R$, we study infinite products of the forms
$$\prod_{n=1}^\infty\Big(R(n)\Big)^{(-1)^{\theta_n}}\quad\text{and}\quad\prod_{n=1}^\infty\Big(R(n)\Big)^{\theta_n}.$$
This generalizes relevant results given by Allouche, Riasat and Shallit in 2019 on infinite products related to the famous Thue-Morse sequence $(t_n)_{n\ge0}$ of the forms
$$\prod_{n=1}^\infty\Big(R(n)\Big)^{(-1)^{t_n}}\quad\text{and}\quad\prod_{n=1}^\infty\Big(R(n)\Big)^{t_n}.$$
\end{abstract}

\maketitle

\section{Introduction}

For any integer $q\ge2$ and $\theta_1,\cdots,\theta_{q-1}\in\{0,1\}$, define the \textit{$(0,\theta_1,\cdots,\theta_{q-1})$-Thue-Morse sequence} $(\theta_n)_{n\ge0}$ to be the unique fixed point of the morphism
$$0\mapsto0\theta_1\cdots\theta_{q-1}$$
$$1\mapsto1\overline{\theta}_1\cdots\overline{\theta}_{q-1}$$
beginning with $\theta_0:=0$, where $\overline{0}:=1$ and $\overline{1}:=0$. Note that the classical Thue-Morse sequence $(t_n)_{n\ge0}$ (see \cite{AS99,AS03}) is exactly the $(0,1)$-Thue-Morse sequence in our terms. Recently Allouche, Riasat and Shallit \cite{ARS19} studied infinite products of the general form $\prod_{n=1}^\infty(R(n))^{(-1)^{t_n}}$ for rational functions $R$, and obtained a class of equalities involving variables in \cite[Theorem 2.2 and Corollary 2.3]{ARS19} together with many concrete equalities in \cite[Corollary2.4]{ARS19}. Besides, they began to study infinite products of the form $\prod_{n=1}^\infty(R(n))^{t_n}$ and obtained three concrete equalities in \cite[Theorem 4.2]{ARS19}. In this paper, we generalize these results by studying infinite products of the forms
$$\prod_{n=1}^\infty\Big(R(n)\Big)^{(-1)^{\theta_n}}\quad\text{and}\quad\prod_{n=1}^\infty\Big(R(n)\Big)^{\theta_n}.$$

Let $\N$, $\N_0$ and $\C$ be the sets of positive integers $1,2,3,\cdots$, non-negative integers $0,1,2,\cdots$ and complex numbers respectively. Moreover, for simplification we denote $\delta_n:=(-1)^{\theta_n}\in\{+1,-1\}$ for all $n\in\N_0$ throughout this paper.

First we have the following convergence theorem, which is a generalization of \cite[Lemma 2.1 and 4.1]{ARS19} (see also \cite[Lemma 1]{R18}) and guarantees the convergence of all the infinite products given in the results in this paper.

\begin{theorem}\label{main-con} Let $q\ge2$ be an integer, $\theta_0=0$, $(\theta_1,\cdots,\theta_{q-1})\in\{0,1\}^{q-1}\setminus\{0^{q-1}\}$, $(\theta_n)_{n\ge0}$ be the $(0,\theta_1,\cdots,\theta_{q-1})$-Thue-Morse sequence and $R\in\C(X)$ be a rational function such that the values $R(n)$ are defined and non-zero for all $n\in\N$. Then:
\begin{itemize}
\item[\emph{(1)}] the infinite product $\prod_{n=1}^\infty (R(n))^{\delta_n}$ converges if and only if the numerator and the denominator of $R$ have the same degree and the same leading coefficient;
\item[\emph{(2)}] the infinite product $\prod_{n=1}^\infty (R(n))^{\theta_n}$ converges if and only if the numerator and the denominator of $R$ have the same degree, the same leading coefficient and the same sum of roots (in $\C$).
\end{itemize}
\end{theorem}

Although Theorem \ref{main-con} is a natural generalization of \cite[Lemma 2.1 and 4.1]{ARS19}, the proof is more intricate and relies on Proposition \ref{partial-sum} as we will see.

In the following Subsection 1.1 and 1.2, we introduce our results on the forms $\prod (R(n))^{\delta_n}$ and $\prod (R(n))^{\theta_n}$ respectively. Then we give some preliminaries in Section 2 and prove all the results in Section 3.

\subsection{Results on the form $\prod (R(n))^{\delta_n}$}

In order to study the infinite product $\prod_{n=1}^\infty(R(n))^{\delta_n}$, by Theorem \ref{main-con} (1), it suffices to study products of the form
$$f(a,b):=\prod_{n=1}^\infty\Big(\frac{n+a}{n+b}\Big)^{\delta_n},$$
where $a,b\in\C\setminus\{-1,-2,-3,\cdots\}$. For the $(0,1)$-Thue-Morse sequence $(t_n)_{n\ge0}$, the special $f(\frac{x}{2},\frac{x+1}{2})=\prod_{n=1}^\infty(\frac{2n+x}{2n+x+1})^{(-1)^{t_n}}$ is used to define new functions in \cite[Theorem 2.2]{ARS19} and \cite[Definition 1]{R18} and further studied (see also \cite[Remark 6.5]{CS13}). For infinite products involving the first $2^m$ terms of $(t_n)_{n\ge0}$, see the equalities (23) and (24) in \cite[Section 6]{CS13}.

As the first main result in this paper, the following theorem generalizes \cite[Theorem 2.2 and Corollary 2.3 (i)]{ARS19} (see also \cite[Lemma 2]{R18} and the equalities (6) and (7) in \cite[Section 4]{R18}).

\begin{theorem}\label{main-1st} Let $q\ge2$ be an integer, $\theta_0=0$, $(\theta_1,\cdots,\theta_{q-1})\in\{0,1\}^{q-1}\setminus\{0^{q-1}\}$ and $(\theta_n)_{n\ge0}$ be the $(0,\theta_1,\cdots,\theta_{q-1})$-Thue-Morse sequence. Then for all $a,b\in\C\setminus\{-1,-2,-3,\cdots\}$, we have
\begin{small}
$$f(a,b)=\Big(\frac{a+1}{b+1}\Big)^{\delta_1}\cdots\Big(\frac{a+q-1}{b+q-1}\Big)^{\delta_{q-1}}f(\frac{a}{q},\frac{b}{q})\Big(f(\frac{a+1}{q},\frac{b+1}{q})\Big)^{\delta_1}\cdots\Big(f(\frac{a+q-1}{q},\frac{b+q-1}{q})\Big)^{\delta_{q-1}},$$
\end{small}
which is equivalent to
$$\prod_{n=1}^\infty\Big(\frac{n+a}{n+b}\cdot\frac{qn+b}{qn+a}\big(\frac{qn+b+1}{qn+a+1}\big)^{\delta_1}\cdots\big(\frac{qn+b+q-1}{qn+a+q-1}\big)^{\delta_{q-1}}\Big)^{\delta_n}=\Big(\frac{a+1}{b+1}\Big)^{\delta_1}\cdots\Big(\frac{a+q-1}{b+q-1}\Big)^{\delta_{q-1}}.$$
\end{theorem}

This theorem implies many neat equalities.

\begin{corollary}\label{cor1-1} Let $q\ge2$ be an integer, $\theta_0=0$, $(\theta_1,\cdots,\theta_{q-1})\in\{0,1\}^{q-1}\setminus\{0^{q-1}\}$ and $(\theta_n)_{n\ge0}$ be the $(0,\theta_1,\cdots,\theta_{q-1})$-Thue-Morse sequence.
\begin{itemize}
\item[\emph{(1)}] For all $a,b\in\C\setminus\{0,-1,-2,\cdots\}$, we have
$$\prod_{n=0}^\infty\Big(\frac{n+a}{n+b}\cdot\frac{qn+b}{qn+a}\big(\frac{qn+b+1}{qn+a+1}\big)^{\delta_1}\cdots\big(\frac{qn+b+q-1}{qn+a+q-1}\big)^{\delta_{q-1}}\Big)^{\delta_n}=1.$$
\item[\emph{(2)}] For all $a\in\C\setminus\{0,-1,-2,\cdots\}$, we have
$$\prod_{n=0}^\infty\Big(\frac{n+a}{n+a+1}\cdot\frac{qn+a+1}{qn+a}\big(\frac{qn+a+2}{qn+a+1}\big)^{\delta_1}\big(\frac{qn+a+3}{qn+a+2}\big)^{\delta_2}\cdots\big(\frac{qn+a+q}{qn+a+q-1}\big)^{\delta_{q-1}}\Big)^{\delta_n}=1$$
and
$$\prod_{n=0}^\infty\Big(\frac{qn+qa}{qn+a}\big(\frac{qn+1}{qn+a+1}\big)^{\delta_1}\big(\frac{qn+2}{qn+a+2}\big)^{\delta_2}\cdots\big(\frac{qn+q-1}{qn+a+q-1}\big)^{\delta_{q-1}}\Big)^{\delta_n}=q.$$
\item[\emph{(3)}] We have
$$\prod_{n=0}^\infty\Big(\frac{qn+q}{qn+1}\big(\frac{qn+1}{qn+2}\big)^{\delta_1}\big(\frac{qn+2}{qn+3}\big)^{\delta_2}\cdots\big(\frac{qn+q-1}{qn+q}\big)^{\delta_{q-1}}\Big)^{\delta_n}=q.$$
\end{itemize}
\end{corollary}

\begin{remark} It seems that Corollary \ref{cor1-1} is weaker than Theorem \ref{main-1st} by the domains of $a$ and $b$. In fact they are equivalent, since Corollary \ref{cor1-1} (1) is the case that $a\neq0$ and $b\neq0$ in Theorem \ref{main-1st}, the second equality in (2) of Corollary \ref{cor1-1} is the case that $a\neq0$ and $b=0$ (the same as $a=0$ and $b\neq0$) in Theorem \ref{main-1st}, and obviously the case that $a=b=0$ in Theorem \ref{main-1st} is trivial.
\end{remark}

Let $q\ge2$ be an integer. For $k=1,2,\cdots,q-1$, define $N_{k,q}(n)$ to be the number of occurrences of the digit $k$ in the base $q$ expansion of the non-negative integer $n$, and let
$$s_q(n):=\sum_{k=1}^{q-1}kN_{k,q}(n)$$
be the sum of digits. It is obtained in \cite[Example 11 and Corollary 5]{AS08} (see also \cite{S06,SG08}) respectively that
\begin{equation}\label{g1}
\prod_{n=0}^\infty\Big(\frac{qn+k}{qn+k+1}\Big)^{(-1)^{N_{k,q}(n)}}=\frac{1}{\sqrt{q}}
\end{equation}
for $k=1,2,\cdots,q-1$, and
\begin{equation}\label{g2}
\prod_{n=0}^\infty\prod_{\substack{0<k<q\\k\text{ odd}}}\Big(\frac{qn+k}{qn+k+1}\Big)^{(-1)^{s_q(n)}}=\frac{1}{\sqrt{q}}.
\end{equation}
For more infinite products related to $(s_q(n))_{n\ge0}$, see for example \cite[Proposition 6 and 7]{L19}. Equalities (\ref{g1}) and (\ref{g2}) are two ways to represent $\frac{1}{\sqrt{q}}$ in the form of infinite products and generalize the well known Woods-Robbins product \cite{R79,W78}
\begin{equation}\label{W-R}
\prod_{n=0}^\infty\Big(\frac{2n+1}{2n+2}\Big)^{(-1)^{t_n}}=\frac{1}{\sqrt{2}}
\end{equation}
where $(t_n)_{n\ge0}$ is the $(0,1)$-Thue-Morse sequence. We give one more such way in the first equality in the following corollary.

\begin{corollary}\label{cor1-2} Let $q\ge2$ be an integer, $k\in\{1,2,\cdots,q-1\}$, $\theta_0=\theta_1=\cdots=\theta_{k-1}=0$, $\theta_k=\theta_{k+1}=\cdots=\theta_{q-1}=1$ and $(\theta_n)_{n\ge0}$ be the $(0,\theta_1,\cdots,\theta_{q-1})$-Thue-Morse sequence. Then
$$\prod_{n=0}^\infty\Big(\frac{qn+k}{qn+q}\Big)^{\delta_n}=\frac{1}{\sqrt{q}}$$
and
$$\prod_{n=0}^\infty\Big(\frac{(n+a)(qn+a+k)^2}{(n+a+1)(qn+a)(qn+a+q)}\Big)^{\delta_n}=1$$
for all $a\in\C\setminus\{0,-1,-2,\cdots\}$.
\end{corollary}

For more generalizations of the Woods-Robbins product (\ref{W-R}), we refer the reader to \cite{AC85,AS89,S85}.

Note that for any integer $q\ge2$, the $(\overbrace{0,\cdots,0}^q)$-Thue-Morse sequence is the trivial $0^\infty$. For $q=2$, the only nontrivial case, related to the $(0,1)$-Thue-Morse sequence, is already studied in \cite{R18} and \cite[Section 2]{ARS19}. In the following three examples, we study nontrivial cases for $q=3$ in detail, related to the $(0,0,1)$, $(0,1,1)$ and $(0,1,0)$-Thue-Morse sequences.

\begin{example}\label{ex1} Let $(\theta_n)_{n\ge0}$ be the $(0,0,1)$-Thue-Morse sequence.
\newline(1) For all $a,b\in\C\setminus\{0,-1,-2,\cdots\}$ we have
$$\prod_{n=0}^\infty\Big(\frac{(n+a)(3n+b)(3n+b+1)(3n+a+2)}{(n+b)(3n+a)(3n+a+1)(3n+b+2)}\Big)^{\delta_n}=1.$$
(2) For all $a\in\C\setminus\{0,-1,-2,\cdots\}$ we have
\begin{equation*}
\begin{aligned}
\textcircled{\tiny{1}}\quad&\prod_{n=0}^\infty\Big(\frac{(n+a)(3n+a+2)^2}{(n+a+1)(3n+a)(3n+a+3)}\Big)^{\delta_n}=1,\\
\textcircled{\tiny{2}}\quad&\prod_{n=0}^\infty\Big(\frac{(3n+1)(3n+3a)(3n+a+2)}{(3n+2)(3n+a)(3n+a+1)}\Big)^{\delta_n}=3,\\
\textcircled{\tiny{3}}\quad&\prod_{n=0}^\infty\Big(\frac{(3n+1)(3n+3a)(3n+a+2)}{(3n+3)(3n+a)(3n+a+1)}\Big)^{\delta_n}=\sqrt{3},\\
\textcircled{\tiny{4}}\quad&\prod_{n=0}^\infty\Big(\frac{(6n+1)(3n+3a)(3n+a+2)}{(6n+5)(3n+a)(3n+a+1)}\Big)^{\delta_n}=1.
\end{aligned}
\end{equation*}
(3) The following concrete equalities hold.
\begin{small}
$$\begin{aligned}
\text{\textcircled{\tiny{1}}}\quad&\prod_{n=0}^\infty\Big(\frac{3n+2}{3n+3}\Big)^{\delta_n}=\frac{1}{\sqrt{3}},\quad&\text{\textcircled{\tiny{2}}}\quad&\prod_{n=0}^\infty\Big(\frac{(6n-3)(6n+3)}{(6n-1)(6n+5)}\Big)^{\delta_n}=1,\\
\text{\textcircled{\tiny{3}}}\quad&\prod_{n=0}^\infty\Big(\frac{(3n+1)(6n+5)}{(3n+2)(6n+1)}\Big)^{\delta_n}=3,\quad&\text{\textcircled{\tiny{4}}}\quad&\prod_{n=0}^\infty\Big(\frac{(3n+1)(6n+5)}{(3n+3)(6n+1)}\Big)^{\delta_n}=\sqrt{3},\\
\text{\textcircled{\tiny{5}}}\quad&\prod_{n=0}^\infty\Big(\frac{(6n+7)^2}{(6n+3)(6n+15)}\Big)^{\delta_n}=1,\quad&\text{\textcircled{\tiny{6}}}\quad&\prod_{n=0}^\infty\Big(\frac{(9n+3)(9n+8)}{(9n+2)(9n+5)}\Big)^{\delta_n}=3,\\
\text{\textcircled{\tiny{7}}}\quad&\prod_{n=0}^\infty\Big(\frac{(18n+3)(18n+17)}{(18n+5)(18n+11)}\Big)^{\delta_n}=1,\quad&\text{\textcircled{\tiny{8}}}\quad&\prod_{n=0}^\infty\Big(\frac{(2n+3)(3n+1)(6n+7)}{(2n+1)(3n+2)(6n+5)}\Big)^{\delta_n}=3,\\
\text{\textcircled{\tiny{9}}}\quad&\prod_{n=0}^\infty\Big(\frac{(n+1)(3n+3)^2}{(n+2)(3n+1)(3n+4)}\Big)^{\delta_n}=1,\quad&\text{\textcircled{\tiny{10}}}\quad&\prod_{n=0}^\infty\Big(\frac{(n+1)(3n+2)(3n+3)}{(n+2)(3n+1)(3n+4)}\Big)^{\delta_n}=\frac{1}{\sqrt{3}},\\
\text{\textcircled{\tiny{11}}}\quad&\prod_{n=0}^\infty\Big(\frac{(n+1)(3n+2)^2}{(n+2)(3n+1)(3n+4)}\Big)^{\delta_n}=\frac{1}{3},\quad&\text{\textcircled{\tiny{12}}}\quad&\prod_{n=0}^\infty\Big(\frac{(3n+2)^3}{(3n+1)(3n+4)(3n+6)}\Big)^{\delta_n}=\frac{1}{3\sqrt{3}},\\
\text{\textcircled{\tiny{13}}}\quad&\prod_{n=0}^\infty\Big(\frac{(n+2)(3n+4)^2}{(n+3)(3n+2)(3n+5)}\Big)^{\delta_n}=1,\quad&\text{\textcircled{\tiny{14}}}\quad&\prod_{n=0}^\infty\Big(\frac{(n+2)(3n+4)^2}{(n+3)(3n+3)(3n+5)}\Big)^{\delta_n}=\frac{1}{\sqrt{3}},\\
\text{\textcircled{\tiny{15}}}\quad&\prod_{n=0}^\infty\Big(\frac{(n+2)(9n+4)(9n+7)}{(n+1)(9n+6)(9n+10)}\Big)^{\delta_n}=1,\quad&\text{\textcircled{\tiny{16}}}\quad&\prod_{n=0}^\infty\Big(\frac{(3n+1)(6n+3)(6n-3)}{(3n+2)(6n+1)(6n-1)}\Big)^{\delta_n}=3.
\end{aligned}$$
\end{small}
\end{example}

\begin{example}\label{ex2} Let $(\theta_n)_{n\ge0}$ be the $(0,1,1)$-Thue-Morse sequence.
\newline(1) For all $a,b\in\C\setminus\{0,-1,-2,\cdots\}$ we have
$$\prod_{n=0}^\infty\Big(\frac{(n+a)(3n+b)(3n+a+1)(3n+a+2)}{(n+b)(3n+a)(3n+b+1)(3n+b+2)}\Big)^{\delta_n}=1.$$
(2) For all $a\in\C\setminus\{0,-1,-2,\cdots\}$ we have
\begin{equation*}
\begin{aligned}
\text{\textcircled{\tiny{1}}}\quad&\prod_{n=0}^\infty\Big(\frac{(n+a)(3n+a+1)^2}{(n+a+1)(3n+a)(3n+a+3)}\Big)^{\delta_n}=1,\\
\text{\textcircled{\tiny{2}}}\quad&\prod_{n=0}^\infty\Big(\frac{(3n+a+1)(3n+a+2)(3n+3a)}{(3n+1)(3n+2)(3n+a)}\Big)^{\delta_n}=3,\\
\text{\textcircled{\tiny{3}}}\quad&\prod_{n=0}^\infty\Big(\frac{(3n+a+1)(3n+a+2)(3n+3a)}{(3n+2)(3n+3)(3n+a)}\Big)^{\delta_n}=\sqrt{3}.
\end{aligned}
\end{equation*}
(3) The following concrete equalities hold.
\begin{small}
$$\begin{aligned}
\text{\textcircled{\tiny{1}}}\quad&\prod_{n=0}^\infty\Big(\frac{3n+1}{3n+3}\Big)^{\delta_n}=\frac{1}{\sqrt{3}},\quad&\text{\textcircled{\tiny{2}}}\quad&\prod_{n=0}^\infty\Big(\frac{(3n+4)(3n+6)}{(3n+2)^2}\Big)^{\delta_n}=\sqrt{3},\\
\text{\textcircled{\tiny{3}}}\quad&\prod_{n=0}^\infty\Big(\frac{(6n+5)^2}{(6n+3)(6n+15)}\Big)^{\delta_n}=1,\quad&\text{\textcircled{\tiny{4}}}\quad&\prod_{n=0}^\infty\Big(\frac{(9n+4)(9n+7)}{(9n+1)(9n+6)}\Big)^{\delta_n}=3,\\
\text{\textcircled{\tiny{5}}}\quad&\prod_{n=0}^\infty\Big(\frac{(9n+5)(9n+8)}{(9n+2)(9n+3)}\Big)^{\delta_n}=3,\quad&\text{\textcircled{\tiny{6}}}\quad&\prod_{n=0}^\infty\Big(\frac{(9n+5)(9n+8)}{(9n+2)(9n+9)}\Big)^{\delta_n}=\sqrt{3},\\
\text{\textcircled{\tiny{7}}}\quad&\prod_{n=0}^\infty\Big(\frac{(24n+7)(24n+13)}{(24n+5)(24n+23)}\Big)^{\delta_n}=1,\quad&\text{\textcircled{\tiny{8}}}\quad&\prod_{n=0}^\infty\Big(\frac{(n+2)(3n+3)^2}{(n+3)(3n+2)(3n+5)}\Big)^{\delta_n}=1,\\
\text{\textcircled{\tiny{9}}}\quad&\prod_{n=0}^\infty\Big(\frac{(n+2)(3n+1)(3n+3)}{(n+3)(3n+2)(3n+5)}\Big)^{\delta_n}=\frac{1}{\sqrt{3}},\quad&\text{\textcircled{\tiny{10}}}\quad&\prod_{n=0}^\infty\Big(\frac{(n+2)(3n+1)^2}{(n+3)(3n+2)(3n+5)}\Big)^{\delta_n}=\frac{1}{3},\\
\text{\textcircled{\tiny{11}}}\quad&\prod_{n=0}^\infty\Big(\frac{(n+3)(3n+4)(3n+5)}{(n+1)(3n+1)(3n+2)}\Big)^{\delta_n}=3,\quad&\text{\textcircled{\tiny{12}}}\quad&\prod_{n=0}^\infty\Big(\frac{(n+3)(3n+4)(3n+5)}{(n+1)(3n+2)(3n+3)}\Big)^{\delta_n}=\sqrt{3},\\
\text{\textcircled{\tiny{13}}}\quad&\prod_{n=0}^\infty\Big(\frac{(3n+4)(3n+5)(3n+9)}{(3n+1)^2(3n+2)}\Big)^{\delta_n}=3\sqrt{3},\quad&\text{\textcircled{\tiny{14}}}\quad&\prod_{n=0}^\infty\Big(\frac{(2n-1)(6n+1)^2}{(2n+1)(6n-1)(6n+5)}\Big)^{\delta_n}=1,\\
\text{\textcircled{\tiny{15}}}\quad&\prod_{n=0}^\infty\Big(\frac{(2n+2)(6n+1)(6n+4)}{(2n+1)(6n+3)(6n+5)}\Big)^{\delta_n}=\frac{1}{\sqrt{3}},\quad&\text{\textcircled{\tiny{16}}}\quad&\prod_{n=0}^\infty\Big(\frac{(2n+3)(6n+5)(6n+7)}{(2n+1)(6n+2)(6n+4)}\Big)^{\delta_n}=3.
\end{aligned}$$
\end{small}
\end{example}

Note that the $(0,1,0)$-Thue-Morse sequence is exactly $01010101\cdots$, which implies $\delta_n:=(-1)^{\theta_n}=(-1)^n$ for all $n\in\N_0$. The next example is deduced from Corollary \ref{cor1-1}, and can also be deduced from Theorem \ref{Gamma prod} and Proposition \ref{Gamma properties}, which are classical results on the Gamma function.

\begin{example}\label{ex3} (1) For all odd $q\ge3$, we have
$$\prod_{n=0}^\infty\Big(\frac{(qn+1)(qn+3)\cdots(qn+q-2)}{(qn+2)(qn+4)\cdots(qn+q-1)}\Big)^{(-1)^n}=\frac{1}{\sqrt{q}}.$$
(2) For all odd $q\ge3$ and all $a\in\C\setminus\{0,-1,-2,\cdots\}$ we have
$$\prod_{n=0}^\infty\Big(\frac{(qn+a)(qn+a+2)(qn+a+4)\cdots(qn+a+q-1)}{(qn+qa)(qn+a+1)(qn+a+3)\cdots(qn+a+q-2)}\Big)^{(-1)^n}=\frac{1}{\sqrt{q}}.$$
(3) The following concrete equalities hold.
\begin{small}
$$\begin{aligned}
\text{\textcircled{\tiny{1}}}\quad&\prod_{n=0}^\infty\Big(\frac{3n+1}{3n+2}\Big)^{(-1)^n}=\frac{1}{\sqrt{3}},\quad&\text{\textcircled{\tiny{2}}}\quad&\prod_{n=0}^\infty\Big(\frac{(n+1)(3n+5)}{(n+3)(3n+4)}\Big)^{(-1)^n}=\frac{1}{\sqrt{3}},\\
\text{\textcircled{\tiny{3}}}\quad&\prod_{n=0}^\infty\Big(\frac{(3n+2)(3n+4)}{(3n+3)(3n+6)}\Big)^{(-1)^n}=\frac{1}{\sqrt{3}},\quad&\text{\textcircled{\tiny{4}}}\quad&\prod_{n=0}^\infty\Big(\frac{(3n+1)(3n+5)}{(3n+6)(3n+9)}\Big)^{(-1)^n}=\frac{1}{3\sqrt{3}},\\
\text{\textcircled{\tiny{5}}}\quad&\prod_{n=0}^\infty\Big(\frac{(9n+2)(9n+8)}{(9n+5)(9n+6)}\Big)^{(-1)^n}=\frac{1}{\sqrt{3}},\quad&\text{\textcircled{\tiny{6}}}\quad&\prod_{n=0}^\infty\Big(\frac{(9n+2)(9n+8)}{(9n+3)(9n+5)}\Big)^{(-1)^n}=1,\\
\text{\textcircled{\tiny{7}}}\quad&\prod_{n=0}^\infty\Big(\frac{(9n+1)(9n+7)}{(9n+3)(9n+4)}\Big)^{(-1)^n}=\frac{1}{\sqrt{3}},\quad&\text{\textcircled{\tiny{8}}}\quad&\prod_{n=0}^\infty\Big(\frac{(9n+1)(9n+7)}{(9n+4)(9n+6)}\Big)^{(-1)^n}=\frac{1}{3},\\
\text{\textcircled{\tiny{9}}}\quad&\prod_{n=0}^\infty\Big(\frac{(9n+2)(9n+11)}{(9n+3)(9n+15)}\Big)^{(-1)^n}=\frac{1}{\sqrt{3}},\quad&\text{\textcircled{\tiny{10}}}\quad&\prod_{n=0}^\infty\Big(\frac{(9n+1)(9n+10)}{(9n+6)(9n+12)}\Big)^{(-1)^n}=\frac{1}{3\sqrt{3}},\\
\text{\textcircled{\tiny{11}}}\quad&\prod_{n=0}^\infty\Big(\frac{(9n+5)(9n+11)}{(9n+8)(9n+15)}\Big)^{(-1)^n}=\frac{1}{\sqrt{3}},\quad&\text{\textcircled{\tiny{12}}}\quad&\prod_{n=0}^\infty\Big(\frac{(9n-1)(9n+8)}{(9n-3)(9n+3)}\Big)^{(-1)^n}=\frac{1}{\sqrt{3}},\\
\text{\textcircled{\tiny{13}}}\quad&\prod_{n=0}^\infty\Big(\frac{(5n+1)(5n+3)}{(5n+2)(5n+4)}\Big)^{(-1)^n}=\frac{1}{\sqrt{5}},\quad&\text{\textcircled{\tiny{14}}}\quad&\prod_{n=0}^\infty\Big(\frac{(10n+1)(10n+9)}{(10n+3)(10n+7)}\Big)^{(-1)^n}=\frac{1}{\sqrt{5}},\\
\text{\textcircled{\tiny{15}}}\quad&\prod_{n=0}^\infty\Big(\frac{(n+1)(5n+3)(5n+7)}{(n+3)(5n+4)(5n+6)}\Big)^{(-1)^n}=\frac{1}{\sqrt{5}},\quad&\text{\textcircled{\tiny{16}}}\quad&\prod_{n=0}^\infty\Big(\frac{(n+1)(5n+2)(5n+7)}{(n+3)(5n+1)(5n+6)}\Big)^{(-1)^n}=1.
\end{aligned}$$
\end{small}
\end{example}

In \cite{H16} Hu studied infinite sums of the form
$$\sum_{n\ge0}\Big((-1)^{a_{w,B}(n)}\sum_{(l,c_l)\in L_{w,B}}c_lf(l(n))\Big)$$
where $a_{w,B}(n)$ denote the number of occurrences of the word $w$ in the base $B$ expansion of the non-negative integer $n$, $f$ is any function that verifies certain convergence conditions, and $L_{w,B}$ is a computable finite set of pairs $(l,c_l)$ where $l$ is a polynomial with integer coefficients of degree $1$ and $c_l$ is an integer. If $f$ is taken to be an appropriate composition of a logarithmic function and a rational function, after exponentiating, some infinite products of the form $\prod_n(R(n))^{(-1)^{a_{w,B}(n)}}$ can be obtained, where $R$ is a rational function depending on the sequence $(a_{w,B}(n))_{n\ge0}$. For instance the above Example \ref{ex3} (3) \textcircled{\tiny{1}} is also obtained in \cite[Section 5]{H16} (see also \cite[Section 4.4]{AS08}).

\subsection{Results on the form $\prod (R(n))^{\theta_n}$}

In order to study the infinite product $\prod_{n=1}^\infty(R(n))^{\theta_n}$, by Theorem \ref{main-con} (2), it suffices to study products of the form
$$\mathfrak{f}(a_1,\cdots,a_d;b_1,\cdots,b_d):=\prod_{n=1}^\infty\Big(\frac{(n+a_1)\cdots(n+a_d)}{(n+b_1)\cdots(n+b_d)}\Big)^{\theta_n}$$
where $d\in\N$ and $a_1,\cdots,a_d,b_1,\cdots,b_d\in\C\setminus\{-1,-2,-3,\cdots\}$ satisfy $a_1+\cdots+a_d=b_1+\cdots+b_d$. As the second main result in this paper, the following theorem (which implies Corollary \ref{TM}) generalizes \cite[Theorem 4.2]{ARS19}.

\begin{theorem}\label{main-2nd} Let $q\ge2$ be an integer, $\theta_0=0$, $(\theta_1,\cdots,\theta_{q-1})\in\{0,1\}^{q-1}\setminus\{0^{q-1}\}$ and $(\theta_n)_{n\ge0}$ be the $(0,\theta_1,\cdots,\theta_{q-1})$-Thue-Morse sequence. Then for all $d\in\N$ and $a_1,\cdots,a_d,b_1,\cdots,b_d\in\C\setminus\{-1,-2,-3,\cdots\}$ satisfying $a_1+\cdots+a_d=b_1+\cdots+b_d$, we have
\begin{small}
$$\mathfrak{f}(a_1,\cdots,a_d;b_1,\cdots,b_d)=\prod_{k=1}^{q-1}\Big(\prod_{i=1}^d\frac{\Gamma(\frac{b_i+k}{q})}{\Gamma(\frac{a_i+k}{q})}\Big)^{\theta_k}\cdot\prod_{k=0}^{q-1}\Big(\mathfrak{f}\big(\frac{a_1+k}{q},\cdots,\frac{a_d+k}{q};\frac{b_1+k}{q},\cdots,\frac{b_d+k}{q}\big)\Big)^{(-1)^{\theta_k}},$$
\end{small}
which is equivalent to
$$\prod_{n=1}^\infty\Big(\prod_{i=1}^d\Big(\frac{n+a_i}{n+b_i}\cdot\prod_{k=0}^{q-1}\big(\frac{qn+b_i+k}{qn+a_i+k}\big)^{(-1)^{\theta_k}}\Big)\Big)^{\theta_n}=\prod_{k=1}^{q-1}\Big(\prod_{i=1}^d\frac{\Gamma(\frac{b_i+k}{q})}{\Gamma(\frac{a_i+k}{q})}\Big)^{\theta_k}.$$
\end{theorem}

This theorem implies a large number of equalities for products of the form $\prod (R(n))^{\theta_n}$ as we will see in the following corollaries, which can also be viewed as special examples.

\begin{corollary}\label{cor2-1} Let $q\ge2$ be an integer, $\theta_0=0$, $(\theta_1,\cdots,\theta_{q-1})\in\{0,1\}^{q-1}\setminus\{0^{q-1}\}$ and $(\theta_n)_{n\ge0}$ be the $(0,\theta_1,\cdots,\theta_{q-1})$-Thue-Morse sequence.
\newline\emph{(1)} For all $a,b,c\in\C$ such that $a,b,a+c,b+c\notin\{-1,-2,-3,\cdots\}$ we have
$$\prod_{n=1}^\infty\Big(\frac{(n+a)(n+b+c)}{(n+b)(n+a+c)}\cdot\prod_{k=0}^{q-1}\Big(\frac{(qn+b+k)(qn+a+c+k)}{(qn+a+k)(qn+b+c+k)}\Big)^{(-1)^{\theta_k}}\Big)^{\theta_n}=\prod_{k=1}^{q-1}\Big(\frac{\Gamma(\frac{b+k}{q})\Gamma(\frac{a+c+k}{q})}{\Gamma(\frac{a+k}{q})\Gamma(\frac{b+c+k}{q})}\Big)^{\theta_k}.$$
\newline\emph{(2)} For all $d\in\N$ and $a_1,\cdots,a_d\in\C\setminus\{-1,-2,-3,\cdots\}$ such that $a_1+\cdots+a_d=0$ we have
$$\prod_{n=1}^\infty\Big(\prod_{i=1}^d\Big(\frac{qn+qa_i}{qn+a_i}\cdot\prod_{k=1}^{q-1}\Big(\frac{qn+k}{qn+a_i+k}\Big)^{(-1)^{\theta_k}}\Big)\Big)^{\theta_n}=\prod_{k=1}^{q-1}\Big(\frac{(\Gamma(\frac{k}{q}))^d}{\Gamma(\frac{a_1+k}{q})\cdots\Gamma(\frac{a_d+k}{q})}\Big)^{\theta_k}.$$
\newline\emph{(3)} For all $a\in\C\setminus\Z$ we have
$$\prod_{n=1}^\infty\Big(\frac{(qn+qa)(qn-qa)}{(qn+a)(qn-a)}\cdot\prod_{k=1}^{q-1}\Big(\frac{(qn+k)^2}{(qn+a+k)(qn-a+k)}\Big)^{(-1)^{\theta_k}}\Big)^{\theta_n}=\prod_{k=1}^{q-1}\Big(\frac{(\Gamma(\frac{k}{q}))^2}{\Gamma(\frac{k+a}{q})\Gamma(\frac{k-a}{q})}\Big)^{\theta_k}.$$
\end{corollary}

In particular for the well known $(0,1)$-Thue-Morse sequence, we have the following corollary, in which (5) \textcircled{\tiny{2}}, \textcircled{\tiny{3}} and \textcircled{\tiny{4}} recover \cite[Theorem 4.2]{ARS19}.

\begin{corollary}\label{TM} Let $(t_n)_{n\ge0}$ be the $(0,1)$-Thue-Morse sequence.
\newline\emph{(1)} For all $d\in\N$ and $a_1,\cdots,a_d,b_1,\cdots,b_d\in\C\setminus\{-1,-2,-3,\cdots\}$ such that $a_1+\cdots+a_d=b_1+\cdots+b_d$ we have
$$\prod_{n=1}^\infty\Big(\prod_{i=1}^d\frac{(n+a_i)(2n+b_i)(2n+a_i+1)}{(n+b_i)(2n+a_i)(2n+b_i+1)}\Big)^{t_n}=\prod_{i=1}^d\frac{\Gamma(\frac{b_i+1}{2})}{\Gamma(\frac{a_i+1}{2})}.$$
\emph{(2)} For all $a,b,c\in\C$ such that $a,b,a+c,b+c\notin\{-1,-2,-3,\cdots\}$ we have
$$\prod_{n=1}^\infty\Big(\frac{(n+a)(n+b+c)(2n+b)(2n+a+1)(2n+a+c)(2n+b+c+1)}{(n+b)(n+a+c)(2n+a)(2n+b+1)(2n+b+c)(2n+a+c+1)}\Big)^{t_n}=\frac{\Gamma(\frac{b+1}{2})\Gamma(\frac{a+c+1}{2})}{\Gamma(\frac{a+1}{2})\Gamma(\frac{b+c+1}{2})}.$$
\begin{itemize}
\item[\emph{(3)} \textcircled{\tiny{1}}] For all $a,b\in\C$ such that $a,b,a+b\notin\{-1,-2,-3,\cdots\}$ we have
$$\prod_{n=1}^\infty\Big(\frac{2(n+a)(n+b)(2n+a+1)(2n+b+1)(2n+a+b)}{(2n+1)(n+a+b)(2n+a)(2n+b)(2n+a+b+1)}\Big)^{t_n}=\frac{\sqrt{\pi}\text{ }\Gamma(\frac{a+b+1}{2})}{\Gamma(\frac{a+1}{2})\Gamma(\frac{b+1}{2})}.$$
\item[\textcircled{\tiny{2}}] For all $a,b\in\C$ such that $a,b,2a+1,a+b\notin\{-1,-2,-3,\cdots\}$ we have
$$\prod_{n=1}^\infty\Big(\frac{(n+a+b)(2n+a+2)(2n+2a+1)(2n+b)(2n+a+b+1)}{(n+2a+1)(2n+a+1)(2n+b+1)(2n+2b)(2n+a+b)}\Big)^{t_n}=\frac{2^a\Gamma(\frac{a+1}{2})\Gamma(\frac{b+1}{2})}{\sqrt{\pi}\text{ }\Gamma(\frac{a+b+1}{2})}.$$
\end{itemize}
\begin{itemize}
\item[\emph{(4)} \textcircled{\tiny{1}}] For all $a\in\C\setminus\{-1,-\frac{3}{2},-2,-\frac{5}{2},\cdots\}$ we have
$$\prod_{n=1}^\infty\Big(\frac{(n+a)(2n+a+2)(2n+2a+1)}{(n+2a+1)(2n+1)(2n+a)}\Big)^{t_n}=2^a.$$
\item[\textcircled{\tiny{2}}] For all $a\in\C\setminus\{-1,-\frac{3}{2},-2,-\frac{5}{2},\cdots\}$ we have
$$\prod_{n=1}^\infty\Big(\frac{(n+1)(n+a+2)(2n+a+3)(2n+2a+1)}{(n+2)(n+2a+1)(2n+3)(2n+a+1)}\Big)^{t_n}=\frac{2^a}{a+1}.$$
\item[\textcircled{\tiny{3}}] For all $a\in\C\setminus\Z$ we have
$$\prod_{n=1}^\infty\Big(\frac{(2n+a+1)(2n-a+1)(2n+2a)(2n-2a)}{(2n+1)^2(2n+a)(2n-a)}\Big)^{t_n}=\cos\frac{\pi a}{2}.$$
\item[\textcircled{\tiny{4}}] For all $a\in\C\setminus(\Z\cup\{\frac{3}{2},\frac{5}{2},\frac{7}{2},\cdots\})$ we have
$$\prod_{n=1}^\infty\Big(\frac{(2n+a+1)(2n-a+1)(2n+2a)(2n-4a+2)}{(2n+1)(2n+a)(2n-a+2)(2n-2a+1)}\Big)^{t_n}=2^a\cos\frac{\pi a}{2}.$$
\item[\textcircled{\tiny{5}}] For all $a\in\C\setminus\{\pm3,\pm5,\pm7,\cdots\}$ we have
$$\prod_{n=1}^\infty\Big(\frac{(2n+a+1)(2n-a+1)(4n+a+3)(4n-a+3)}{(2n+2)^2(4n+a+1)(4n-a+1)}\Big)^{t_n}=\frac{\sqrt{\pi}}{\Gamma(\frac{3+a}{4})\Gamma(\frac{3-a}{4})}.$$
\item[\textcircled{\tiny{6}}] For all $d\in\N$ we have
$$\prod_{n=1}^\infty\Big(\frac{(n+1)(2n+d)(2n+2)^{2d-1}}{(n+d)(2n+d+1)(2n+1)^{2d-1}}\Big)^{t_n}=\pi^{\frac{d-1}{2}}\Gamma(\frac{d+1}{2}).$$
\end{itemize}
\emph{(5)} The following concrete equalities hold.
\begin{footnotesize}
$$\begin{aligned}
\text{\textcircled{\tiny{1}}}\quad&\prod_{n=0}^\infty\Big(\frac{(2n+1)(4n-1)}{(2n-1)(4n+3)}\Big)^{t_n}=\sqrt{2},\quad&\text{\textcircled{\tiny{2}}}\quad&\prod_{n=0}^\infty\Big(\frac{(2n+1)(4n+3)}{(2n+2)(4n+1)}\Big)^{t_n}=\frac{\Gamma(\frac{1}{4})}{\sqrt{2}\pi^{\frac{3}{4}}},\\
\text{\textcircled{\tiny{3}}}\quad&\prod_{n=0}^\infty\Big(\frac{(n+1)(4n+5)}{(n+2)(4n+1)}\Big)^{t_n}=\sqrt{2},\quad&\text{\textcircled{\tiny{4}}}\quad&\prod_{n=0}^\infty\Big(\frac{(8n+1)(8n+7)}{(8n+3)(8n+5)}\Big)^{t_n}=\sqrt{2\sqrt{2}-2},\\
\text{\textcircled{\tiny{5}}}\quad&\prod_{n=0}^\infty\Big(\frac{(n+1)(2n+3)^2}{(n+3)(2n+1)^2}\Big)^{t_n}=2,\quad&\text{\textcircled{\tiny{6}}}\quad&\prod_{n=0}^\infty\Big(\frac{(3n+2)^2(6n+5)}{(3n+3)^2(6n+1)}\Big)^{t_n}=\frac{\sqrt{3}\text{ }\Gamma(\frac{1}{3})\Gamma(\frac{1}{6})}{4\pi^{\frac{3}{2}}},\\
\text{\textcircled{\tiny{7}}}\quad&\prod_{n=0}^\infty\Big(\frac{(n+2)^2(2n+5)}{(n+1)(n+5)(2n+1)}\Big)^{t_n}=4,\quad&\text{\textcircled{\tiny{8}}}\quad&\prod_{n=0}^\infty\Big(\frac{(2n+1)^2(4n-1)}{(2n-1)(2n+2)(4n+1)}\Big)^{t_n}=\frac{\Gamma(\frac{1}{4})}{\pi^{\frac{3}{4}}},\\
\text{\textcircled{\tiny{9}}}\quad&\prod_{n=0}^\infty\Big(\frac{(2n+3)^2(4n-1)}{(2n-1)(2n+6)(4n+1)}\Big)^{t_n}=\frac{2\Gamma(\frac{1}{4})}{\pi^{\frac{3}{4}}},\quad&\text{\textcircled{\tiny{10}}}\quad&\prod_{n=0}^\infty\Big(\frac{(2n-1)(4n+3)^2}{(2n+2)(4n+1)(4n-1)}\Big)^{t_n}=\frac{\Gamma(\frac{1}{4})}{2\pi^{\frac{3}{4}}},\\
\text{\textcircled{\tiny{11}}}\quad&\prod_{n=0}^\infty\Big(\frac{(3n-1)^2(6n+3)}{(3n+2)(3n-2)(6n-1)}\Big)^{t_n}=2^{\frac{2}{3}},\quad&\text{\textcircled{\tiny{12}}}\quad&\prod_{n=0}^\infty\Big(\frac{(4n+2)^2(8n-1)}{(4n-1)(4n+1)(8n+7)}\Big)^{t_n}=2^{\frac{1}{4}},\\
\text{\textcircled{\tiny{13}}}\quad&\prod_{n=0}^\infty\Big(\frac{(n+1)(2n+7)(4n+9)}{(n+4)(2n+3)(4n+5)}\Big)^{t_n}=\frac{4\sqrt{2}}{5},\quad&\text{\textcircled{\tiny{14}}}\quad&\prod_{n=0}^\infty\Big(\frac{(n+1)(3n+7)(6n+5)}{(n+2)(3n+2)(6n+9)}\Big)^{t_n}=3\cdot2^{-\frac{5}{3}},\\
\text{\textcircled{\tiny{15}}}\quad&\prod_{n=0}^\infty\Big(\frac{(3n+1)(6n-1)(6n+3)}{(3n-1)(6n+1)(6n+5)}\Big)^{t_n}=2^{\frac{1}{3}},\quad&\text{\textcircled{\tiny{16}}}\quad&\prod_{n=0}^\infty\Big(\frac{(5n+4)(10n+1)(10n+5)}{(5n+2)(10n+3)(10n+7)}\Big)^{t_n}=\frac{\sqrt{5}-1}{2^{\frac{2}{5}}}.\\
\end{aligned}$$
\end{footnotesize}
\end{corollary}

\section{Preliminaries}

First we need the following concept.

\begin{definition}[\cite{AS08,U99}]
Let $q\ge2$ be an integer. A sequence $u=(u_n)_{n\ge0}\in\C^{\N_0}$ is called \textit{strongly $q$-multiplicative} if $u_0=1$ and
$$u_{nq+k}=u_nu_k$$
for all $k\in\{0,1,\cdots,q-1\}$ and $n\in\N_0$.
\end{definition}

The following theorem is a classical result on the Gamma function $\Gamma$ (see for examples \cite[Theorem 1.1]{CS13} and \cite[Section 12.13]{WW96}).

\begin{theorem}\label{Gamma prod} Let $d\in\N$ and $a_1,a_2,\cdots,a_d,b_1,b_2,\cdots,b_d\in\C\setminus\{0,-1,-2,\cdots\}$. If $a_1+a_2+\cdots+a_d=b_1+b_2+\cdots+b_d$, then
$$\prod_{n=0}^\infty\frac{(n+a_1)(n+a_2)\cdots(n+a_d)}{(n+b_1)(n+b_2)\cdots(n+b_d)}=\frac{\Gamma(b_1)\Gamma(b_2)\cdots\Gamma(b_d)}{\Gamma(a_1)\Gamma(a_2)\cdots\Gamma(a_d)}.$$
\end{theorem}

Besides, we need the properties on the Gamma function gathered in the following proposition.

\begin{proposition}[\cite{A64,B17,W}]\label{Gamma properties} \emph{(1)} For all $n\in\N$ and $z\in\C\setminus\{0,-\frac{1}{n},-\frac{2}{n},-\frac{3}{n},-\frac{4}{n},\cdots\}$ we have
$$\Gamma(z)\Gamma(z+\frac{1}{n})\Gamma(z+\frac{2}{n})\cdots\Gamma(z+\frac{n-1}{n})=(2\pi)^{\frac{n-1}{2}}n^{\frac{1}{2}-nz}\Gamma(nz).$$
\emph{(2)} For all $z\in\C\setminus\{0,-1,-2,\cdots\}$ we have
$$\Gamma(z+1)=z\Gamma(z)$$
and
$$\Gamma(\frac{z}{2})\Gamma(\frac{z+1}{2})=2^{1-z}\sqrt{\pi}\text{ }\Gamma(z).$$
\emph{(3)} For all $z\in\C\setminus\Z$ we have
$$\Gamma(z)\Gamma(1-z)=\frac{\pi}{\sin\pi z}.$$
\emph{(4)} We have
$$\Gamma(1)=\Gamma(2)=1,\quad\Gamma(\frac{1}{2})=\sqrt{\pi}\quad\text{and}\quad\Gamma(\frac{3}{2})=\frac{\sqrt{\pi}}{2}.$$
\end{proposition}

\section{Proofs of the results}

Let $q\ge2$ be an integer, $\theta_0=0$, $\theta_1,\cdots,\theta_{q-1}\in\{0,1\}$ and $(\theta_n)_{n\ge0}$ be the $(0,\theta_1,\cdots,\theta_{q-1})$-Thue-Morse sequence. Recall that $(\delta_n)_{n\ge0}$ is defined by $\delta_n=(-1)^{\theta_n}$ for all $n\in\N_0$. At the same time $(\delta_n)_{n\ge0}$ can be view as the unique fixed point of the morphism
\begin{equation}\label{+-morphism}
\begin{aligned}
+1\mapsto(+1)(+\delta_1)\cdots(+\delta_{q-1})\\
-1\mapsto(-1)(-\delta_1)\cdots(-\delta_{q-1})
\end{aligned}
\end{equation}
beginning with $\delta_0=+1$. Define the sequence of partial sums of $(\delta_n)_{n\ge0}$ by
$$\Delta_0:=0\quad\text{and}\quad\Delta_n:=\delta_0+\delta_1+\cdots+\delta_{n-1}\quad\text{for all }n\ge1.$$
Note that $(\Delta_n)_{n\ge0}$ depends on the choice of $(\delta_1,\cdots,\delta_{q-1})\in\{+1,-1\}^{q-1}$. Before proving Theorem \ref{main-con}, we need the following proposition on $(\Delta_n)_{n\ge0}$, which is itself valuable.

\begin{proposition}\label{partial-sum} Let $q\ge2$ be an integer.
\begin{itemize}
\item[\emph{(1)}] For all $k,s\in\N_0$ and $t\in\{0,1,\cdots,q^k-1\}$ we have
$$\delta_{sq^k+t}=\delta_s\delta_t\quad\text{and}\quad\Delta_{sq^k+t}=\Delta_s\Delta_{q^k}+\delta_s\Delta_t.$$
\item[\emph{(2)}] With the convention $0^0:=1$, for all $k\in\N_0$ we have
$$\Delta_{q^k}=\Delta_q^k,$$
$$\max\Big\{|\Delta_{q^k}|:(\delta_1,\cdots,\delta_{q-1})\in\{+1,-1\}^{q-1}\setminus\{(+1)^{q-1}\}\Big\}=(q-2)^k,$$
$$\max\Big\{|\Delta_n|:0\le n\le q^k,(\delta_1,\cdots,\delta_{q-1})\in\{+1,-1\}^{q-1}\setminus\{(+1)^{q-1}\}\Big\}=1+(q-2)+\cdots+(q-2)^k.$$
\item[\emph{(3)}] If $(\delta_1,\cdots,\delta_{q-1})\neq(+1)^{q-1}$, then for all $n$ large enough we have
$$|\Delta_n|\le n^{\log_q(q-1)}.$$
\end{itemize}
\end{proposition}
\begin{proof} (1) \textcircled{\tiny{1}} Prove $\delta_{sq^k+t}=\delta_s\delta_t$ for all $k,s\in\N_0$ and $t\in\{0,1,\cdots,q^k-1\}$.
\begin{itemize}
\item[i)] Prove that $(\delta_n)_{n\ge0}$ is strongly $q$-multiplicative, i.e.,
$$\delta_{sq+t}=\delta_s\delta_t\quad\text{for all }s\in\N_0\text{ and }t\in\{0,1,\cdots,q-1\}.$$
Denote the morphism (\ref{+-morphism}) by $\psi$. Then by $\psi((\delta_0,\delta_1,\delta_2,\cdots))=(\delta_0,\delta_1,\delta_2,\cdots)$ we get $\psi(\delta_s)=(\delta_{sq},\delta_{sq+1},\cdots,\delta_{(s+1)q-1})$ for all $s\in\N_0$. It follows from $\psi(+1)=(+1,+\delta_1,\cdots,+\delta_{q-1})$ and $\psi(-1)=(-1,-\delta_1,\cdots,-\delta_{q-1})$ that $\delta_{sq+t}=\delta_s\delta_t$ for all $t\in\{0,\cdots,q-1\}$.
\item[ii)] Let $k\in\N$, $s\in\N_0$ and $t\in\{0,\cdots,q^k-1\}$. Then there exist $l\in\N_0$ and $s_l,\cdots,s_1,s_0$, $t_{k-1},\cdots,t_1,t_0$ $\in\{0,1,\cdots,q-1\}$ such that
$$s=s_lq^l+\cdots+s_1q+s_0\quad\text{and}\quad t=t_{k-1}q^{k-1}+\cdots+t_1q+t_0.$$
By i) and \cite[Proposition 1]{AS08} we get
$$\delta_{sq^k+t}=\delta_{s_l}\cdots\delta_{s_1}\delta_{s_0}\delta_{t_{k-1}}\cdots\delta_{t_1}\delta_{t_0},$$
$$\delta_s=\delta_{s_l}\cdots\delta_{s_1}\delta_{s_0}\quad\text{and}\quad\delta_t=\delta_{t_{k-1}}\cdots\delta_{t_1}\delta_{t_0}.$$
Thus $\delta_{sq^k+t}=\delta_s\delta_t$.
\end{itemize}
\textcircled{\tiny{2}} Prove $\Delta_{sq^k+t}=\Delta_s\Delta_{q^k}+\delta_s\Delta_t$ for all $k,s\in\N_0$ and $t\in\{0,1,\cdots,q^k-1\}$. In fact, we have
\begin{eqnarray*}
\Delta_{sq^k+t}&=&(\delta_0+\delta_1+\cdots+\delta_{q^k-1})+(\delta_{q^k}+\delta_{q^k+1}+\cdots+\delta_{q^k+(q^k-1)})\\
& &+\cdots+(\delta_{(s-1)q^k}+\delta_{(s-1)q^k+1}+\cdots+\delta_{(s-1)q^k+(q^k-1)})\\
& &+(\delta_{sq^k}+\delta_{sq^k+1}+\cdots+\delta_{sq^k+t-1})\\
&=&\delta_0(\delta_0+\delta_1+\cdots+\delta_{q^k-1})+\delta_1(\delta_0+\delta_1+\cdots+\delta_{q^k-1})\\
& &+\cdots+\delta_{s-1}(\delta_0+\delta_1+\cdots+\delta_{q^k-1})\\
& &+\delta_s(\delta_0+\delta_1+\cdots+\delta_{t-1})\\
&=&\Delta_s\Delta_{q^k}+\delta_s\Delta_t
\end{eqnarray*}
where the second equality follows from \textcircled{\tiny{1}}.
\newline(2) \textcircled{\tiny{1}} We have $\Delta_{q^k}=\Delta_q^k$ for all $k\in\N_0$ since (1) \textcircled{\tiny{2}} implies $\Delta_{q\cdot q^l}=\Delta_q\Delta_{q^l}$ for all $l\in\N_0$.
\newline\textcircled{\tiny{2}} For all $k\in\N_0$, the fact
$$\max\Big\{|\Delta_{q^k}|:(\delta_1,\cdots,\delta_{q-1})\in\{+1,-1\}^{q-1}\setminus\{(+1)^{q-1}\}\Big\}=(q-2)^k$$
follows from \textcircled{\tiny{1}} and
$$\max\Big\{|\Delta_q|:(\delta_1,\cdots,\delta_{q-1})\in\{+1,-1\}^{q-1}\setminus\{(+1)^{q-1}\}\Big\}=q-2.$$
\textcircled{\tiny{3}} In order to prove the last equality in statement (2), since the case $k=0$ is trivial and \textcircled{\tiny{2}} implies $|\Delta_{q^k}|\le1+(q-2)+(q-2)^2+\cdots+(q-2)^k$, it suffices to verify that for all $k\in\N$, we have
\newline\leftline{$\max\Big\{|\Delta_n|:0\le n\le q^k-1,(\delta_1,\cdots,\delta_{q-1})\in\{+1,-1\}^{q-1}\setminus\{(+1)^{q-1}\}\Big\}$}
\rightline{$=1+(q-2)+(q-2)^2+\cdots+(q-2)^k.$}
\begin{itemize}
\item[\boxed{\le}] (By induction on $k$) For $k=1$, obviously we have $|\Delta_0|,|\Delta_1|,\cdots,|\Delta_{q-1}|\le q-1$. Suppose that for some $k\in\N$ and all $l\in\{0,1,\cdots,k\}$, we have already had
$$|\Delta_0|,|\Delta_1|,\cdots,|\Delta_{q^l-1}|\le1+(q-2)+(q-2)^2+\cdots+(q-2)^l.$$
Let $n\in\{0,1,\cdots,q^{k+1}-1\}$. It suffices to prove
\begin{equation}\label{to prove}
|\Delta_n|\le1+(q-2)+(q-2)^2+\cdots+(q-2)^{k+1}.
\end{equation}
If $n\le q^k-1$, this follows immediately from the inductive hypothesis. We only need to consider $q^k\le n\le q^{k+1}-1$ in the following. Let $s\in\{1,\cdots,q-1\}$ and $t\in\{0,1,\cdots,q^k-1\}$ such that $n=sq^k+t$. By (1) \textcircled{\tiny{2}} we get
$$\Delta_n=\Delta_s\Delta_{q^k}+\delta_s\Delta_t.$$
If $s\le q-2$, then
\begin{eqnarray*}
|\Delta_n|&\le&|\Delta_s|\cdot|\Delta_{q^k}|+|\Delta_t|\\
&\le&s(q-2)^k+(1+(q-2)+(q-2)^2+\cdots+(q-2)^k)\\
&\le&1+(q-2)+(q-2)^2+\cdots+(q-2)^{k+1}
\end{eqnarray*}
where the second inequality follows from \textcircled{\tiny{2}} and the inductive hypothesis. In the following we only need to consider $s=q-1$. It means that
$$\Delta_n=\Delta_{q-1}\Delta_{q^k}+\delta_{q-1}\Delta_t.$$
If there exists $p\in\{0,1,\cdots,q-2\}$ such that $\delta_p=-1$, then $|\Delta_{q-1}|\le q-3$ and
\begin{eqnarray*}
|\Delta_n|&\le&(q-3)|\Delta_{q^k}|+|\Delta_t|\\
&\le&(q-3)(q-2)^k+(1+(q-2)+(q-2)^2+\cdots+(q-2)^k)\\
&\le&1+(q-2)+(q-2)^2+\cdots+(q-2)^{k+1}
\end{eqnarray*}
where the second inequality follows from \textcircled{\tiny{2}} and the inductive hypothesis. Thus it suffices to consider $\delta_0=\delta_1=\cdots=\delta_{q-2}=+1$ in the following. By $(\delta_1,\cdots,\delta_{q-1})\neq(+1)^{q-1}$ we get $\delta_{q-1}=-1$. It follows from $\Delta_{q-1}=q-1$ and $\Delta_{q^k}=\Delta_q^k=(q-2)^k$ that
$$\Delta_n=(q-1)(q-2)^k-\Delta_t.$$
Thus proving (\ref{to prove}) is equivalent to verifying
$$-1-(q-2)-\cdots-(q-2)^{k-1}\le\Delta_t\le1+(q-2)+\cdots+(q-2)^{k+1}+(q-1)(q-2)^k.$$
Since the second inequality follows immediately from the inductive hypothesis, we only need to prove the first inequality. Let $u\in\{0,1,\cdots,q-1\}$ and $v\in\{0,1,\cdots,q^{k-1}-1\}$ such that $t=uq^{k-1}+v$. By (1) \textcircled{\tiny{2}} we get
$$\Delta_t=\Delta_u\Delta_{q^{k-1}}+\delta_u\Delta_v.$$
Since $\delta_0=\delta_1=\cdots=\delta_{q-2}=+1$, $\delta_{q-1}=-1$ and $0\le u\le q-1$ imply $\Delta_u=u\ge0$, $\Delta_q=q-2$ and $\Delta_{q^{k-1}}=\Delta_q^{k-1}=(q-2)^{k-1}\ge0$, by $\delta_u\in\{+1,-1\}$ we get
$$\Delta_t\ge-|\Delta_v|\ge-1-(q-2)-\cdots-(q-2)^{k-1}$$
where the last inequality follows from the inductive hypothesis.
\item[\boxed{\ge}] Let $\delta_1=\delta_2=\cdots=\delta_{q-2}=+1$ and $\delta_{q-1}=-1$. It suffices to prove that for all $k\in\N$ we have
\begin{equation}\label{to prove example}
\Delta_{q^k-q^{k-1}-\cdots-q-1}=(q-2)^k+\cdots+(q-2)^2+(q-2)+1.
\end{equation}
(By induction) For $k=1$ we have $\Delta_{q-1}=q-1$. Suppose that (\ref{to prove example}) is true for some $k\in\N$. Then for $k+1$, we have
\begin{eqnarray*}
\Delta_{q^{k+1}-q^k-q^{k-1}-\cdots-q-1}&=&\Delta_{(q-2)q^k+(q^k-q^{k-1}-\cdots-q-1)}\\
&=&\Delta_{q-2}\Delta_{q^k}+\delta_{q-2}\Delta_{q^k-q^{k-1}-\cdots-q-1}\\
&=&(q-2)\Delta_q^k+(q-2)^k+\cdots+(q-2)^2+(q-2)+1\\
&=&(q-2)^{k+1}+(q-2)^k+\cdots+(q-2)^2+(q-2)+1
\end{eqnarray*}
where the second equality follows from (1) \textcircled{\tiny{2}} and the third equality follows from \textcircled{\tiny{1}} and the inductive hypothesis.
\end{itemize}
(3) For $n\in\N$ large enough, there exists $k\in\N$ large enough such that $q^k+1\le n\le q^{k+1}$. By (2) \textcircled{\tiny{3}} we get
$$|\Delta_n|\le1+(q-2)+\cdots+(q-2)^{k+1}\le(q-1)^k=(q^k)^{\log_q(q-1)}\le n^{\log_q(q-1)}$$
where the second inequality can be verified straightforwardly for $k$ large enough.
\end{proof}

\begin{proof}[Proof of Theorem \ref{main-con}] Since (2) follows in the same way as in the proof of \cite[Lemma 4.1]{ARS19} by applying (1), we only need to prove (1) in the following.
\newline\boxed{\Rightarrow} Suppose that $\prod_{n=1}^\infty (R(n))^{\delta_n}$ converges. Then $(R(n))^{\delta_n}\to1$ as $n\to\infty$. Since $\delta_n\in\{+1,-1\}$ for all $n\in\N$, we get $R(n)\to1$ as $n\to\infty$. Thus the numerator and the denominator of $R$ have the same degree and the same leading coefficient.
\newline\boxed{\Leftarrow} Suppose that the numerator and the denominator of $R$ have the same leading coefficient and the same degree. Decompose them into factors of degree $1$. To prove that $\prod_{n=1}^\infty (R(n))^{\delta_n}$ converges, it suffices to show that $\prod_{n=1}^\infty(\frac{n+a}{n+b})^{\delta_n}$ converges for all $a,b\in\C$ satisfying $n+a\neq0$ and $n+b\neq0$ for all $n\in\N$ (that is, $a,b\in\C\setminus\{-1,-2,-3,\cdots\}$). Since $(\frac{n+a}{n+b})^{\delta_n}\to1$ as $n\to\infty$, we only need to prove that
$$\prod_{n=1}^\infty\Big(\big(\frac{qn+a}{qn+b}\big)^{\delta_{qn}}\big(\frac{qn+1+a}{qn+1+b}\big)^{\delta_{qn+1}}\cdots\big(\frac{qn+q-1+a}{qn+q-1+b}\big)^{\delta_{qn+q-1}}\Big)$$
converges. Since Proposition \ref{partial-sum} (1) implies $\delta_{qn}=\delta_n\delta_0$, $\delta_{qn+1}=\delta_n\delta_1$, $\cdots$, $\delta_{qn+q-1}=\delta_n\delta_{q-1}$, it suffices to show that
$$\prod_{n=1}^\infty\big(r(n)\big)^{\delta_n}$$
converges, where
$$r(n):=\big(\frac{qn+a}{qn+b}\big)^{\delta_0}\big(\frac{qn+1+a}{qn+1+b}\big)^{\delta_1}\cdots\big(\frac{qn+q-1+a}{qn+q-1+b}\big)^{\delta_{q-1}}.$$
This is equivalent to show that
\begin{equation}\label{sum}
\sum_{n=1}^\infty\delta_n\ln r(n)
\end{equation}
converges. Since there exist $c_0,c_1,\cdots,c_{q-1},d_0,d_1,\cdots,d_{q-1}\in\C$ such that
\begin{small}
$$r(n)=\frac{q^qn^q+c_{q-1}n^{q-1}+\cdots+c_1n+c_0}{q^qn^q+d_{q-1}n^{q-1}+\cdots+d_1n+d_0}=1+\frac{(c_{q-1}-d_{q-1})n^{q-1}+\cdots+(c_1-d_1)n+(c_0-d_0)}{q^qn^q+d_{q-1}n^{q-1}+\cdots+d_1n+d_0},$$
\end{small}
we get
$$\ln r(n)-\frac{c_{q-1}-d_{q-1}}{q^qn}=\cO(\frac{1}{n^2}),$$
which implies that
$$\sum_{n=1}^\infty\delta_n\big(\ln r(n)-\frac{c_{q-1}-d_{q-1}}{q^qn}\big)$$
converges absolutely. In order to prove that (\ref{sum}) converges, we only need to show that
$$\sum_{n=1}^\infty\frac{\delta_n}{n}$$
converges. Enlightened by partial summation (see for example the equality (6.5) in \cite{BMS20} related to the Thue-Morse sequence), we consider the following \textcircled{\tiny{1}} and \textcircled{\tiny{2}}, which complete the proof.
\begin{itemize}
\item[\textcircled{\tiny{1}}] Prove that
$$\sum_{n=1}^\infty\frac{\delta_1+\cdots+\delta_n}{n(n+1)}$$
converges. In fact, since Proposition \ref{partial-sum} (3) implies
$$\frac{|\Delta_n|}{n^2}\le\frac{1}{n^{2-\log_q(q-1)}}\quad\text{for all }n\text{ large enough},$$
where $2-\log_q(q-1)>1$, it follows that $\sum_{n=1}^\infty\frac{\Delta_n}{n^2}$ converges absolutely. So does $\sum_{n=1}^\infty\frac{\Delta_n}{n(n+1)}$. Thus we only need to check that $\sum_{n=1}^\infty(\frac{\delta_1+\cdots+\delta_n}{n(n+1)}-\frac{\Delta_n}{n(n+1)})$ converges. This follows immediately from $|\delta_1+\cdots+\delta_n-\Delta_n|=|\delta_n-\delta_0|\le2$.
\item[\textcircled{\tiny{2}}] Prove that
$$\sum_{n=1}^\infty(\frac{\delta_n}{n}-\frac{\delta_1+\cdots+\delta_n}{n(n+1)})$$
converges to $0$. In fact, for all $N\in\N$ we have
\begin{eqnarray*}
\sum_{n=1}^N\frac{\delta_1+\cdots+\delta_n}{n(n+1)}&=&\sum_{n=1}^N(\delta_1+\cdots+\delta_n)(\frac{1}{n}-\frac{1}{n+1})\\
&=&\delta_1\sum_{n=1}^N(\frac{1}{n}-\frac{1}{n+1})+\delta_2\sum_{n=2}^N(\frac{1}{n}-\frac{1}{n+1})+\cdots+\delta_N\sum_{n=N}^N(\frac{1}{n}-\frac{1}{n+1})\\
&=&\delta_1(1-\frac{1}{N+1})+\delta_2(\frac{1}{2}-\frac{1}{N+1})+\cdots+\delta_N(\frac{1}{N}-\frac{1}{N+1})\\
&=&\sum_{n=1}^N\frac{\delta_n}{n}-\frac{\delta_1+\delta_2+\cdots+\delta_N}{N+1},
\end{eqnarray*}
which implies
$$\sum_{n=1}^N(\frac{\delta_n}{n}-\frac{\delta_1+\cdots+\delta_n}{n(n+1)})=\frac{\delta_1+\delta_2+\cdots+\delta_N}{N+1}=\frac{\Delta_{N+1}-1}{N+1}.$$
Since Proposition \ref{partial-sum} (3) implies
$$\frac{|\Delta_{N+1}|}{N+1}\le\frac{1}{(N+1)^{1-\log_q(q-1)}}\quad\text{for all }N\text{ large enough},$$
where $1-\log_q(q-1)>0$, as $N\to\infty$ we get $\frac{\Delta_{N+1}}{N+1}\to0$ and then $\sum_{n=1}^N(\frac{\delta_n}{n}-\frac{\delta_1+\cdots+\delta_n}{n(n+1)})\to0$.
\end{itemize}
\end{proof}

\begin{proof}[Proof of Theorem \ref{main-1st}] Since Proposition \ref{partial-sum} (1) implies $\delta_{qn}=\delta_n\delta_0,\delta_{qn+1}=\delta_n\delta_1,\cdots,\delta_{qn+q-1}=\delta_n\delta_{q-1}$ for all $n\in\N_0$, we get $f(a,b)$
\begin{small}
$$\begin{aligned}
&=\prod_{n=1}^\infty\Big(\frac{qn+a}{qn+b}\Big)^{\delta_{qn}}\prod_{n=0}^\infty\Big(\frac{qn+1+a}{qn+1+b}\Big)^{\delta_{qn+1}}\cdots\prod_{n=0}^\infty\Big(\frac{qn+q-1+a}{qn+q-1+b}\Big)^{\delta_{qn+q-1}}\\
&=\prod_{n=1}^\infty\Big(\frac{qn+a}{qn+b}\Big)^{\delta_n\delta_0}\prod_{n=0}^\infty\Big(\frac{qn+a+1}{qn+b+1}\Big)^{\delta_n\delta_1}\cdots\prod_{n=0}^\infty\Big(\frac{qn+a+q-1}{qn+b+q-1}\Big)^{\delta_n\delta_{q-1}}\\
&=\Big(\frac{a+1}{b+1}\Big)^{\delta_0\delta_1}\cdots\Big(\frac{a+q-1}{b+q-1}\Big)^{\delta_0\delta_{q-1}}\prod_{n=1}^\infty\Big(\frac{n+\frac{a}{q}}{n+\frac{b}{q}}\Big)^{\delta_n\delta_0}\prod_{n=1}^\infty\Big(\frac{n+\frac{a+1}{q}}{n+\frac{b+1}{q}}\Big)^{\delta_n\delta_1}\cdots\prod_{n=1}^\infty\Big(\frac{n+\frac{a+q-1}{q}}{n+\frac{b+q-1}{q}}\Big)^{\delta_n\delta_{q-1}}\\
&=\Big(\frac{a+1}{b+1}\Big)^{\delta_1}\cdots\Big(\frac{a+q-1}{b+q-1}\Big)^{\delta_{q-1}}f(\frac{a}{q},\frac{b}{q})\Big(f(\frac{a+1}{q},\frac{b+1}{q})\Big)^{\delta_1}\cdots\Big(f(\frac{a+q-1}{q},\frac{b+q-1}{q})\Big)^{\delta_{q-1}}.
\end{aligned}$$
\end{small}
\end{proof}

\begin{proof}[Proof of Corollary \ref{cor1-1}] (1) follows from Theorem \ref{main-1st} after multiplying by the factor corresponding to $n=0$. The first equality in (2) follows from taking $b=a+1$ in (1). The second equality in (2) follows from taking $b=0$ in Theorem \ref{main-1st} and then multiplying the factor corresponding to $n=0$. We should note that it does not follow from taking $b=0$ in (1). Finally (3) follows immediately from taking $a=1$ in the second equality in (2).
\end{proof}

\begin{proof}[Proof of Corollary \ref{cor1-2}] These two equalities follow from Corollary \ref{cor1-1} (3) and the first equality in (2) of Corollary \ref{cor1-1} respectively.
\end{proof}

\begin{proof}[Proof of Example \ref{ex1}] (1) follows from Corollary \ref{cor1-1} (1).
\begin{itemize}
\item[(2)]\textcircled{\tiny{1}} and \textcircled{\tiny{2}} follow from Corollary \ref{cor1-1} (2).
\item[]\textcircled{\tiny{3}} follows from \textcircled{\tiny{2}} and the fact that the first equality in Corollary \ref{cor1-2} implies
\begin{equation}\label{above}
\prod_{n=0}^\infty\Big(\frac{3n+2}{3n+3}\Big)^{\delta_n}=\frac{1}{\sqrt{3}}.
\end{equation}
\item[]\textcircled{\tiny{4}} follows from taking $b=\frac{1}{2}$ in (1).
\item[(3)] \textcircled{\tiny{1}} is the above equality (\ref{above}).
\item[]\textcircled{\tiny{2}}, \textcircled{\tiny{5}}, \textcircled{\tiny{9}} and \textcircled{\tiny{13}} follow from taking $a=-\frac{1}{2},\frac{3}{2},1$ and $2$ respectively in (2) \textcircled{\tiny{1}}.
\item[]\textcircled{\tiny{3}}, \textcircled{\tiny{6}}, \textcircled{\tiny{8}} and \textcircled{\tiny{16}} follow from taking $a=\frac{1}{2},\frac{2}{3},\frac{3}{2}$ and $-\frac{1}{2}$ respectively in (2) \textcircled{\tiny{2}}.
\item[]\textcircled{\tiny{4}} follows from multiplying \textcircled{\tiny{3}} and \textcircled{\tiny{1}}.
\item[]\textcircled{\tiny{7}} follows from taking $a=\frac{5}{6}$ in (2) \textcircled{\tiny{4}}.
\item[]\textcircled{\tiny{10}}, \textcircled{\tiny{11}}, \textcircled{\tiny{12}} and \textcircled{\tiny{14}} follow respectively from \textcircled{\tiny{9}}, \textcircled{\tiny{10}}, \textcircled{\tiny{11}} and \textcircled{\tiny{13}} by applying \textcircled{\tiny{1}}.
\item[]\textcircled{\tiny{15}} follows from taking $a=2$ and $b=\frac{4}{3}$ in (1).
\end{itemize}
\end{proof}

\begin{proof}[Proof of Example \ref{ex2}] (1) follows from Corollary \ref{cor1-1} (1).
\begin{itemize}
\item[(2)]\textcircled{\tiny{1}} and \textcircled{\tiny{2}} follow from Corollary \ref{cor1-1} (2).
\item[]\textcircled{\tiny{3}} follows from \textcircled{\tiny{2}} and the fact that the first equality in Corollary \ref{cor1-2} implies
$$\prod_{n=0}^\infty\Big(\frac{3n+1}{3n+3}\Big)^{\delta_n}=\frac{1}{\sqrt{3}}.$$
\item[(3)]\textcircled{\tiny{1}} is the above equality.
\item[]\textcircled{\tiny{2}} follows from taking $a=2$ in (2) \textcircled{\tiny{3}}.
\item[]\textcircled{\tiny{3}} and \textcircled{\tiny{8}} follow from taking $a=\frac{3}{2}$ and $2$ respectively in (2) \textcircled{\tiny{1}}.
\item[]\textcircled{\tiny{4}}, \textcircled{\tiny{5}}, \textcircled{\tiny{11}} and \textcircled{\tiny{16}} follow from taking $a=\frac{1}{3},\frac{2}{3},3$ and $\frac{3}{2}$ respectively in (2) \textcircled{\tiny{2}}.
\item[]\textcircled{\tiny{6}}, \textcircled{\tiny{9}} and \textcircled{\tiny{10}} follow respectively from \textcircled{\tiny{5}}, \textcircled{\tiny{8}} and \textcircled{\tiny{9}} by applying \textcircled{\tiny{1}}.
\item[]\textcircled{\tiny{7}} follows from taking $a=\frac{5}{8}$ and $b=\frac{7}{8}$ in (1).
\item[]\textcircled{\tiny{12}} and \textcircled{\tiny{13}} follow respectively from multiplying and dividing \textcircled{\tiny{11}} by \textcircled{\tiny{1}}.
\item[]\textcircled{\tiny{14}} follows from combining the results of taking $a=\frac{1}{2}$ and $-\frac{1}{2}$ in (2) \textcircled{\tiny{2}}.
\item[]\textcircled{\tiny{15}} follows from taking $a=1,b=\frac{1}{2}$ in (1) and then multiplying by \textcircled{\tiny{1}}.
\end{itemize}
\end{proof}

\begin{proof}[Proof of Example \ref{ex3}] For odd $q\ge3$, let $\theta_1=\theta_3=\cdots=\theta_{q-2}=1$ and $\theta_2=\theta_4=\cdots=\theta_{q-1}=0$. Then the $(0,\theta_1,\cdots,\theta_{q-1})$-Thue-Morse sequence $(\theta_n)_{n\ge0}$ is exactly $(01)^\infty$. It follows that $\delta_n:=(-1)^{\theta_n}=(-1)^n$ for all $n\ge0$.
\newline(1) By the second equality in Corollary \ref{cor1-1} (2) we get
\begin{footnotesize}
\begin{equation}\label{q}
\prod_{n=0}^\infty\Big(\frac{(qn+qa)(qn+a+1)(qn+2)(qn+a+3)(qn+4)\cdots(qn+a+q-2)(qn+q-1)}{(qn+a)(qn+1)(qn+a+2)(qn+3)(qn+a+4)\cdots(qn+q-2)(qn+a+q-1)}\Big)^{(-1)^n}=q
\end{equation}
\end{footnotesize}
for all $a\in\C\setminus\{0,-1,-2,\cdots\}$. Then we conclude (1) by taking $a=1$ in (\ref{q}).
\newline(2) follows from (\ref{q}) and (1).
\newline(3) Note that for all $q\in\N$ and $a\in\C\setminus\{0,-1,-2,\cdots\}$ we have
\begin{equation}\label{qa}
\prod_{n=0}^\infty\Big(\frac{(qn+a)(qn+a+q)}{(qn+qa)(qn+qa+q)}\Big)^{(-1)^n}=\frac{1}{q}
\end{equation}
since the left hand side is
$$\lim_{k\to\infty}\frac{a}{qa}\cdot\frac{a+q}{qa+q}\cdot\Big(\frac{a+q}{qa+q}\cdot\frac{a+2q}{qa+2q}\Big)^{-1}\cdot\frac{a+2q}{qa+2q}\cdot\frac{a+3q}{qa+3q}\cdots\Big(\frac{a+kq}{qa+kq}\cdot\frac{a+(k+1)q}{qa+(k+1)q}\Big)^{(-1)^k}$$
$$=\lim_{k\to\infty}\frac{a}{qa}\cdot\Big(\frac{a+(k+1)q}{qa+(k+1)q}\Big)^{(-1)^k}=\frac{1}{q}.$$
We prove the concrete equalities in the following.
\newline\textcircled{\tiny{1}} and \textcircled{\tiny{13}} follow from taking $q=3$ and $5$ respectively in (1).
\newline\textcircled{\tiny{2}}, \textcircled{\tiny{3}}, \textcircled{\tiny{5}} and \textcircled{\tiny{7}} follow from taking $q=3$, and then $a=3,2,\frac{2}{3}$ and $\frac{1}{3}$ respectively in (2).
\newline\textcircled{\tiny{4}}, \textcircled{\tiny{9}}, \textcircled{\tiny{10}} and \textcircled{\tiny{12}} are deduced by applying \textcircled{\tiny{1}} noting that (\ref{qa}) with $q=3$ and then $a=2,\frac{2}{3},\frac{1}{3}$ and $-\frac{1}{3}$ give respectively
    $$\prod_{n=0}^\infty\Big(\frac{(3n+2)(3n+5)}{(3n+6)(3n+9)}\Big)^{(-1)^n}=\frac{1}{3},\quad\prod_{n=0}^\infty\Big(\frac{(9n+2)(9n+11)}{(9n+6)(9n+15)}\Big)^{(-1)^n}=\frac{1}{3},$$
    $$\prod_{n=0}^\infty\Big(\frac{(9n+1)(9n+10)}{(9n+3)(9n+12)}\Big)^{(-1)^n}=\frac{1}{3}\quad\text{and}\quad\prod_{n=0}^\infty\Big(\frac{(9n-1)(9n+8)}{(9n-3)(9n+6)}\Big)^{(-1)^n}=\frac{1}{3}.$$
\newline\textcircled{\tiny{6}}, \textcircled{\tiny{8}}, \textcircled{\tiny{11}} and \textcircled{\tiny{16}} follow respectively from dividing \textcircled{\tiny{5}} by \textcircled{\tiny{1}}, multiplying \textcircled{\tiny{7}} by \textcircled{\tiny{1}}, dividing \textcircled{\tiny{9}} by \textcircled{\tiny{6}} and dividing \textcircled{\tiny{15}} by \textcircled{\tiny{13}}.
\newline\textcircled{\tiny{14}} and \textcircled{\tiny{15}} follow from taking $q=5$, and then $a=\frac{1}{2}$ and $3$ respectively in (2).
\end{proof}

Before proving Theorem \ref{main-2nd}, we need the following proposition.

\begin{proposition}\label{theta qn+k} Let $q\ge2$ be an integer, $\theta_0=0$, $(\theta_1,\cdots,\theta_{q-1})\in\{0,1\}^{q-1}\setminus\{0^{q-1}\}$ and $(\theta_n)_{n\ge0}$ be the $(0,\theta_1,\cdots,\theta_{q-1})$-Thue-Morse sequence. Then for all $n\in\N_0$ and $k\in\{0,1,\cdots,q-1\}$ we have
$$\theta_{nq+k}=\theta_n(-1)^{\theta_k}+\theta_k.$$
\end{proposition}
\begin{proof} Let $h$ denote the morphism
$$0\mapsto0\theta_1\cdots\theta_{q-1}$$
$$1\mapsto1\overline{\theta}_1\cdots\overline{\theta}_{q-1}$$
where $\overline{0}:=1$ and $\overline{1}:=0$. By $h(\theta_0\theta_1\theta_2\cdots)=\theta_0\theta_1\theta_2\cdots$ we get
$$h(\theta_n)=\theta_{nq}\theta_{nq+1}\cdots\theta_{nq+q-1}$$
for all $n\in\N_0$. It follows from $h(0)=\theta_0\theta_1\cdots\theta_{q-1}$ and $h(1)=\overline{\theta}_0\overline{\theta}_1\cdots\overline{\theta}_{q-1}$ that
$$\theta_{nq+k}=\left\{\begin{array}{ll}
\theta_k & \text{if } \theta_n=0\\
\overline{\theta}_k & \text{if } \theta_n=1
\end{array}\right.\quad=\theta_n(-1)^{\theta_k}+\theta_k\quad\text{for all }k\in\{0,1,\cdots,q-1\}.$$
\end{proof}

\begin{proof}[Proof of Theorem \ref{main-2nd}] We have $\mathfrak{f}(a_1,\cdots,a_d;b_1,\cdots,b_d)$
\begin{eqnarray*}
&=&\prod_{n=1}^\infty\Big(\prod_{i=1}^d\frac{n+a_i}{n+b_i}\Big)^{\theta_n}\\
&=&\prod_{k=1}^{q-1}\Big(\prod_{i=1}^d\frac{k+a_i}{k+b_i}\Big)^{\theta_k}\cdot\prod_{n=1}^\infty\prod_{k=0}^{q-1}\Big(\prod_{i=1}^d\frac{nq+k+a_i}{nq+k+b_i}\Big)^{\theta_{nq+k}}\\
&\overset{(\star)}{=}&\prod_{k=1}^{q-1}\Big(\prod_{i=1}^d\frac{a_i+k}{b_i+k}\Big)^{\theta_k}\cdot\prod_{n=1}^\infty\prod_{k=0}^{q-1}\Big(\prod_{i=1}^d\frac{qn+a_i+k}{qn+b_i+k}\Big)^{\theta_n(-1)^{\theta_k}+\theta_k}\\
&=&\prod_{k=1}^{q-1}\Big(\prod_{i=1}^d\frac{a_i+k}{b_i+k}\Big)^{\theta_k}\cdot\prod_{n=1}^\infty\prod_{k=0}^{q-1}\Big(\prod_{i=1}^d\frac{qn+a_i+k}{qn+b_i+k}\Big)^{\theta_k}\cdot\prod_{n=1}^\infty\prod_{k=0}^{q-1}\Big(\prod_{i=1}^d\frac{qn+a_i+k}{qn+b_i+k}\Big)^{\theta_n(-1)^{\theta_k}}\\
&\overset{(\star\star)}{=}&\prod_{k=1}^{q-1}\Big(\prod_{i=1}^d\frac{a_i+k}{b_i+k}\Big)^{\theta_k}\cdot\prod_{n=1}^\infty\prod_{k=1}^{q-1}\Big(\prod_{i=1}^d\frac{qn+a_i+k}{qn+b_i+k}\Big)^{\theta_k}\cdot\prod_{k=0}^{q-1}\prod_{n=1}^\infty\Big(\prod_{i=1}^d\frac{qn+a_i+k}{qn+b_i+k}\Big)^{\theta_n(-1)^{\theta_k}}\\
&=&\prod_{n=0}^\infty\prod_{k=1}^{q-1}\Big(\prod_{i=1}^d\frac{qn+a_i+k}{qn+b_i+k}\Big)^{\theta_k}\cdot\prod_{k=0}^{q-1}\Big(\prod_{n=1}^\infty\Big(\prod_{i=1}^d\frac{qn+a_i+k}{qn+b_i+k}\Big)^{\theta_n}\Big)^{(-1)^{\theta_k}}\\
&=&\prod_{k=1}^{q-1}\Big(\prod_{n=0}^\infty\prod_{i=1}^d\frac{n+\frac{a_i+k}{q}}{n+\frac{b_i+k}{q}}\Big)^{\theta_k}\cdot\prod_{k=0}^{q-1}\Big(\prod_{n=1}^\infty\Big(\prod_{i=1}^d\frac{n+\frac{a_i+k}{q}}{n+\frac{b_i+k}{q}}\Big)^{\theta_n}\Big)^{(-1)^{\theta_k}}\\
&\overset{(\star\star\star)}{=}&\prod_{k=1}^{q-1}\Big(\prod_{i=1}^d\frac{\Gamma(\frac{b_i+k}{q})}{\Gamma(\frac{a_i+k}{q})}\Big)^{\theta_k}\cdot\prod_{k=0}^{q-1}\Big(\mathfrak{f}\big(\frac{a_1+k}{q},\cdots,\frac{a_d+k}{q};\frac{b_1+k}{q},\cdots,\frac{b_d+k}{q}\big)\Big)^{(-1)^{\theta_k}},
\end{eqnarray*}
where ($\star$), ($\star\star$) and ($\star\star\star$) follow from Proposition \ref{theta qn+k}, $\theta_0=0$ and Theorem \ref{Gamma prod} respectively.
\end{proof}

\begin{proof}[Proof of Corollary \ref{cor2-1}] (1) follows from taking $d=2$, $a_1=a$, $a_2=b+c$, $b_1=b$ and $b_2=a+c$ in Theorem \ref{main-2nd}.
\newline(2) follows from taking $b_1=\cdots=b_d=0$ in Theorem \ref{main-2nd}.
\newline(3) follows from taking $d=2$, $a_1=a$ and $a_2=-a$ in (2).
\end{proof}

\begin{proof}[Proof of Corollary \ref{TM}] In the following proof, for calculations related to the Gamma function, we use Proposition \ref{Gamma properties} frequently without citations. (1) and (2) follow from Theorem \ref{main-2nd} and Corollary \ref{cor2-1} (1) respectively.
\begin{itemize}
\item[(3)]\textcircled{\tiny{1}} follows from taking $b=0$ in (2) and then replacing all $c$ by $b$.
\item[]\textcircled{\tiny{2}} follows from taking $c=a-1$ in (2) and then replacing all $a$ by $a+1$.
\item[(4)]\textcircled{\tiny{1}} follows from multiplying (3) \textcircled{\tiny{1}} and \textcircled{\tiny{2}}.
\item[]\textcircled{\tiny{2}} follows from taking $b=2$ in (3) \textcircled{\tiny{2}}.
\item[]\textcircled{\tiny{3}} and \textcircled{\tiny{4}} follow from taking $b=-a$ and $1-2a$ respectively in (3) \textcircled{\tiny{1}}.
\item[]\textcircled{\tiny{5}} follows from taking $d=2,a_1=\frac{1+a}{2},a_2=\frac{1-a}{2},b_1=0$ and $b_2=1$ in (1).
\item[]\textcircled{\tiny{6}} follows from taking $a_1=\cdots=a_d=1,b_1=d$ and $b_2=\cdots=b_d=0$ in (1).
\item[(5)]\textcircled{\tiny{1}} follows from taking $a=\frac{1}{2}$ in (4) \textcircled{\tiny{3}}.
\item[]\textcircled{\tiny{2}} and \textcircled{\tiny{6}} follow from taking $a=0$ and $\frac{1}{3}$ respectively in (4) \textcircled{\tiny{5}}.
\item[]\textcircled{\tiny{3}}, \textcircled{\tiny{5}}, \textcircled{\tiny{7}}, \textcircled{\tiny{11}} and \textcircled{\tiny{12}} follow from taking $a=\frac{1}{2}$, $1$, $2$, $-\frac{2}{3}$ and $-\frac{1}{4}$ respectively in (4) \textcircled{\tiny{1}}.
\item[]\textcircled{\tiny{4}}, \textcircled{\tiny{15}} and \textcircled{\tiny{16}} follow from taking $a=\frac{1}{4}$, $\frac{2}{3}$ and $\frac{2}{5}$ respectively in (4) \textcircled{\tiny{4}}.
\item[]\textcircled{\tiny{8}}, \textcircled{\tiny{9}} and \textcircled{\tiny{10}} follow respectively from multiplying \textcircled{\tiny{1}} by \textcircled{\tiny{2}}, multiplying \textcircled{\tiny{5}} by \textcircled{\tiny{8}} and dividing \textcircled{\tiny{2}} by \textcircled{\tiny{1}}.
\item[]\textcircled{\tiny{13}} and \textcircled{\tiny{14}} follow from taking $a=\frac{3}{2}$ and $\frac{1}{3}$ respectively in (4) \textcircled{\tiny{2}}.
\end{itemize}
\end{proof}

\begin{ack}
The author thanks Prof. Jean-Paul Allouche for his advices, and thanks the Oversea Study Program of Guangzhou Elite Project (GEP) for financial support (JY201815). While the author was preparing this paper, he learned that Dr. Shuo Li was working on infinite products related to $\phi$-Thue-Morse sequence, which is another generalization of the classical Thue-Morse sequence, and in particular for the classical one, some new equalities are obtained.
\end{ack}

\end{document}